\documentclass[11.05pt,a4paper]{amsart}
\usepackage{amssymb,amsmath}
\usepackage{latexsym}
\usepackage{amsthm,amsfonts,amssymb,mathrsfs}
\usepackage{rotating}
\usepackage[leqno]{amsmath}
\usepackage{xspace}
\usepackage[all]{xy}
\usepackage{longtable}

\headsep=1cm \numberwithin{equation}{section}
\newtheorem{theorem}{Theorem}[section]
\newtheorem{lemma}[theorem]{Lemma}
\newtheorem{proposition}[theorem]{Proposition}
\newtheorem{corollary}[theorem]{Corollary}

\theoremstyle{definition}

\theoremstyle{remark}

\newtheorem{conjecture}[theorem]{Conjecture}
\newtheorem{acknowledgement}{Acknowledgement}

\newcommand{\Ass}{\operatorname{Ass}}
\newcommand{\im}{\operatorname{im}}
\newcommand{\Ker}{\operatorname{Ker}}

\newcommand{\id}{\operatorname{id}}

\newcommand{\Gid}{\operatorname{Gid}}
\newcommand{\pd}{\operatorname{pd}}
\newcommand{\Gpd}{\operatorname{Gpd}}

\newcommand{\V}{\operatorname{V}}
\newcommand{\G}{\operatorname{G}}
\newcommand{\E}{\operatorname{E}}

\newcommand{\Ext}{\operatorname{Ext}}
\newcommand{\Supp}{\operatorname{Supp}}

\newcommand{\Hom}{\operatorname{Hom}}
\newcommand{\RHom}{\operatorname{RHom}}

\newcommand{\Coker}{\operatorname{Coker}}

\newcommand{\Max}{\operatorname{Max}}

\newcommand{\lo}{\longrightarrow}
\newcommand{\fm}{\frak{m}}
\newcommand{\fp}{\frak{p}}

\newcommand{\fn}{\frak{n}}

\newenvironment{prf}[1][Proof]{\begin{proof}[\bf #1]}{\end{proof}}

\begin{document}

\author[M. Nikkhah Babaei and K. Divaani-Aazar]{Massoumeh Nikkhah Babaei and Kamran Divaani-Aazar}
\title[Gorenstein injective envelopes ...]
{Gorenstein injective envelopes of Artinian modules}

\address{M. Nikkhah Babaei, Department of Mathematics, Alzahra University, Vanak, Post Code 19834,
Tehran, Iran.}
\email{massnikkhah@yahoo.com}

\address{K. Divaani-Aazar, Department of Mathematics, Alzahra University, Vanak, Post Code
19834, Tehran, Iran-and-School of Mathematics, Institute for Research in Fundamental Sciences
(IPM), P.O. Box 19395-5746, Tehran, Iran.}
\email{kdivaani@ipm.ir}

\subjclass[2010]{13D05; 13E10.}

\keywords{Gorenstein injective dimension; Gorenstein injective envelopes; Gorenstein injective
modules; Gorenstein projective covers; Gorenstein projective modules; injective covers; Henselian
rings.\\
The second author was supported by a grant from IPM (No. 91130212).}

\begin{abstract} Let $R$ be a commutative Noetherian ring and $A$ an Artinian $R$-module. We
prove that if $A$ has finite Gorenstein injective dimension, then $A$ possesses a Gorenstein
injective envelope which is special and Artinian. This, in particular, yields that over a
Gorenstein ring any Artinian module possesses a Gorenstein injective envelope which is special
and Artinian.
\end{abstract}

\maketitle

\section{Introduction and Prerequisites}

Throughout this paper, $R$ is a commutative Noetherian ring with an identity. It is a standard fact that
any $R$-module $M$ possesses an injective envelope $\E_R(M)$. Contrary to that it is not known yet that
any $R$-module has a Gorenstein injective envelope. For an $R$-module $M$, let $\G_R(M)$ denote the
Gorenstein injective envelope of $M$ if it exists. It is known that if $A$ is an Artinian $R$-module, then
$\E_R(A)$ is also Artinian. There is much evidence supporting Holm's claim ``Every result in classical
homological algebra should have a counterpart in Gorenstein homological algebra''; see \cite{H}.
Hence, the following question arises: Does any Artinian $R$-module $A$ admit an Artinian Gorenstein
injective envelope? In this paper, we are concerned with this question. Belshoff and Enochs
\cite[Theorem 1]{BE} have given a partial positive answer to this question. They proved that if $R$ is a
Gorenstein local ring with the residue field $k$, then $\G_R(k)$ is Artinian.  We show that if $A$ is an
Artinian $R$-module with finite Gorenstein injective dimension, then $\G_R(A)$ exists, it is Artinian and
$\id_R(\G_R(A)/A)<\infty$. In particular, we conclude that if $R$ is Gorenstein, then any Artinian
$R$-module admits an Artinian Gorenstein injective envelope. The following are two of the remaining
open questions in the theory of Gorenstein dimensions; see \cite[Questions 3.17 and 3.19]{CFoH}:\\
Is every direct limit of Gorenstein injective $R$-modules Gorenstein injective?\\
Is $M_{\fp}$ Gorenstein injective over $R_{\fp}$ for all Gorenstein injective $R$-modules $M$ and
all prime ideals $\fp$ of $R$?\\
Assuming a positive answer to the first question, we provide an affirmative answer to the second
question.

Below, we recall some definitions and notations.

An $R$-module $M$ is said to be \emph{Gorenstein injective} if there exists an exact sequence
$$I^{\bullet}: \cdots \lo I_1\lo I_0\lo I^0\lo I^1\lo \cdots$$ of injective $R$-modules such
that $M\cong \im(I_0\lo I^0)$ and that $\Hom_R(I,I^{\bullet})$ is exact for all injective
$R$-modules $I$. The notion of \emph{Gorenstein injective dimension} of an $R$-module $M$,
$\Gid M$, is defined as the infimum of the length of right resolutions of $M$ which are
consisting of Gorenstein injective $R$-modules.

Let $\mathscr{X}$ be a class of $R$-modules and $M$ an $R$-module. An $R$-homomorphism
$\phi:M\lo X$ where $X\in \mathscr{X}$ is called a $\mathscr{X}$-\emph{preenvelope}
of $M$ if for any  $X'\in \mathscr{X}$, the induced $R$-homomorphism $\Hom_R(X,X')\lo \Hom_R(M,X')$
is surjective. If, moreover, any $R$-homomorphism $f:X\lo X$ such that $f \phi=\phi$ is an
automorphism, then $\phi$ is called a $\mathscr{X}$-\emph{envelope} of $M$. From the definition, it
becomes clear that if $\mathscr{X}$-envelopes exist, then they are unique up to isomorphism.  A 
$\mathscr{X}$-(pre)envelope $\phi:M\lo X$ is called \emph{special} if $\Ext_R^i(\Coker\phi,X')=0$
for all $X'\in \mathscr{X}$ and all $i>0$.

If $\mathscr{X}$ is the class of Gorenstein injective $R$-modules, then a $\mathscr{X}$-(pre)envelope
is called \emph{Gorenstein injective (pre)envelope}. By \cite[Corollary 2.7]{EL}, over a commutative
Noetherian ring, any module has a Gorenstein injective preenvelope. In fact by \cite[Theorem 7.12]{K} 
any such a module has a special Gorenstein injective preenvelope. But, it is not known whether any
module over an arbitrary commutative Noetherian ring has a Gorenstein injective envelope or not. Since
any injective $R$-module is Gorenstein injective, it follows that any Gorenstein injective preenvelope
is a monomorphism. If $M$ admits a Gorenstein injective envelope, then without loss of generality,
we may assume that $M$ has a Gorenstein injective envelope $i:M\hookrightarrow \G_R(M)$, where
$M$ is a submodule of the Gorenstein injective $R$-module $\G_R(M)$ and $i$ is the inclusion
map. For simplicity, sometimes we call $\G_R(M)$ the Gorenstein injective envelope of $M$.  If
$\phi:M\lo E$ is a Gorenstein injective (pre)envelope such that $\Coker \phi$ has finite injective
dimension, then using \cite[Lemma 2.2]{CFrH}, we can deduce that $\phi:M\lo E$ is a special Gorenstein
injective (pre)envelope of $M$.

An $R$-module $M$ is said to be \emph{Gorenstein projective} if there exists an exact sequence
$$P_{\bullet}:\cdots \lo P_1\lo P_0\lo P^0\lo P^1\lo \cdots $$ of projective $R$-modules such
that $M\cong \im (P_0\lo P^0)$ and that $\Hom_R(P_{\bullet},P)$ is exact for all projective $R$-
modules $P$. The \emph{Gorenstein projective dimension} of an $R$-module $M$, $\Gpd_RM$, is defined
as the infimum of the length of left resolutions of $M$ which are consisting of Gorenstein projective
$R$-modules.

Let $\mathscr{X}$ be a class of $R$-modules and $M$ an $R$-module. An $R$-homomorphism $\phi:X\lo M$ where
$X\in \mathscr{X}$ is called a $\mathscr{X}$-\emph{precover} of $M$ if for any $X'\in \mathscr{X}$, the
induced $R$-homomorphism $\Hom_R(X',X)\lo \Hom_R(X',M)$ is surjective. If, moreover, any $R$-homomorphism
$f:X\lo X$ such that $\phi f=\phi$ is an automorphism, then $\phi$ is called a $\mathscr{X}$-\emph{cover}
of $M$. Obviously, if $\mathscr{X}$-covers exist, then they are unique up to isomorphism. A
$\mathscr{X}$-(pre)cover $\phi:X\lo M$ is called \emph{special} if $\Ext_R^i(X',\Ker \phi)=0$ for all
$X'\in \mathscr{X}$ and all $i>0$. If $\mathscr{X}$ is the class of Gorenstein projective $R$-modules,
then a $\mathscr{X}$-(pre)cover is called \emph{Gorenstein projective (pre)cover}. As any projective
$R$-module is Gorenstein projective, it follows that any Gorenstein projective precover is an epimorphism.
If $\mathscr{X}$ is the class of injective $R$-modules, then a $\mathscr{X}$-(pre)cover is called
\emph{injective (pre)cover}. By \cite[Theorem 5.4.1 and Corollary 7.2.3]{EJ}, over a commutative Noetherian
ring, every $R$-module admits a special injective cover.

\section{The Results}

Theorem \ref{210} is our main result. To prove it, we need to the following five lemmas. We start by the
following  special case of \cite[Corollary 2.5]{T2} that can be easily deduced from \cite[Corollary 2.13]{H}
and \cite[Remark 2.6]{T1}.

\begin{lemma}\label{21} Let $(R,\fm)$ be a Henselian local ring and $N$ a finitely generated $R$-module
of finite Gorenstein projective dimension. Then $N$ possesses a special Gorenstein projective cover
$\phi:V\lo N$ such that $V$ is finitely generated and $\Ker \phi$ has finite projective dimension.
\end{lemma}

In what follows, the following will serve us as a technical tool.

\begin{lemma}\label{23} Let $f:R\lo T$ be a homomorphism of Noetherian commutative rings and $M$ a
$T$-module. Assume that any Gorenstein injective $T$-module is also a Gorenstein injective $R$-module via $f$.
\begin{enumerate}
\item[i)] $\Gid_RM\leq \Gid_TM$ and equality holds if any $T$-module which is Gorenstein injective
as an $R$-module is also Gorenstein injective as a $T$-module.
\item[ii)]  Assume $f$ is flat. Let $\phi:M\lo E$ be a Gorenstein injective preenvelope
of $M$ as a $T$-module such that $\id_T(\Coker \phi)<\infty$. Then $\phi:M\lo E$ is also a Gorenstein
injective preenvelope of $M$ as an $R$-module and $\id_R(\Coker \phi)<\infty$.
\end{enumerate}
\end{lemma}

\begin{prf} i) Our assumption implies that any Gorenstein injective resolution of $M$ as a $T$-module is
also a Gorenstein injective resolution of $M$ as an $R$-module. Hence,  $\Gid_RM\leq \Gid_TM$.

Now, assume that any $T$-module which is Gorenstein injective as an $R$-module is also Gorenstein
injective as a $T$-module. We show that $\Gid_TM\leq \Gid_RM$. We may assume that $n:=\Gid_RM$ is finite.
Let $$E^{\bullet}: 0\lo E^0\lo E^1\lo \cdots \lo E^i\lo \cdots$$ be a Gorenstein injective resolution of
$M$ as a $T$-module. Since $E^{\bullet}$ is also a Gorenstein injective resolution of $M$ as an $R$-module,
\cite[Theorem 3.3]{CFrH} yields that $\Ker(E^n\lo E^{n+1})$ is a Gorenstein injective $R$-module.
Then, by the assumption $\Ker(E^n\lo E^{n+1})$ is also Gorenstein injective as a $T$-module. Thus, we
deduce that $$0\lo E^0\lo E^1\lo \cdots \lo E^{n-1}\lo \Ker(E^n\lo E^{n+1})\lo 0$$ is a Gorenstein injective
resolution of $M$ as a $T$-module, and so $\Gid_TM\leq \Gid_RM$.

ii) By the assumption $E$ is also a Gorenstein injective $R$-module. As any injective $T$-module is also
injective as an $R$-module, we deduce that $\id_R(\Coker \phi)<\infty$. Hence, \cite[Lemma 2.2]{CFrH} implies
that $\Ext_R^i(\Coker \phi,G)=0$ for all Gorenstein injective $R$-modules $G$ and all $i>0$. Therefore, for
any Gorenstein injective $R$-module $G$, from the exact sequence $$0\lo M\overset{\phi}\lo E\lo \Coker \phi
\lo 0,$$ we obtain the exact sequence $$\Hom_R(E,G)\lo \Hom_R(M,G)\lo \Ext_R^1(\Coker \phi,G)=0.$$ This
completes the argument.
\end{prf}

Let $(R,\fm)$ be a local ring. To be able to apply Lemma \ref{23} ii) to  the natural flat ring homomorphism
$\theta:R\lo \hat{R}$, we have to prove the following:

\begin{lemma}\label{24} Let $f:(R,\fm)\lo (T,\fn)$ be a faithfully flat ring homomorphism of local rings.
Then every Gorenstein injective $T$-module is also a Gorenstein injective $R$-module via $f$.
\end{lemma}

\begin{prf} Let $M$ be a Gorenstein injective $T$-module. As $\pd_{T}(T\otimes_RT)<\infty$,
\cite[Lemma 2.2]{CFrH} yields that $\Ext^i_{T}(T\otimes_RT,M)=0$ for all $i>0$. On the other
hand, using adjointness, one can check that $\Ext^i_R(T,M)\cong \Ext^i_{T}(T\otimes_RT,M)$
for all $i\geq 0$. So, $\RHom_R(T,M)\simeq \Hom_R(T,M)$. Since, by \cite[Ascent table
\textbf{I}, (d)]{CH}, $\Hom_R(T,M)$ is a Gorenstein injective $T$-module, we deduce that
$\Gid_{T}(\RHom_R(T,M))=0$. Now, \cite[Theorem 1.7]{CW} implies that $\Gid_RM=\Gid_{T}
(\RHom_R(T,M)),$  and so $M$ is a Gorenstein injective $R$-module.
\end{prf}

Let $A$ be an Artinian module over a local ring $(R,\fm)$.  Sharp \cite{Sh} has shown that $A$ has a
natural structure as a module over $\hat{R}$. Let $\theta:R\lo \hat{R}$ denote the natural ring homomorphism.
The $\hat{R}$-module structure of $A$ is such that for any element $r\in R$ the multiplication by $r$ on $A$
has the same effect as the multiplication of $\theta(r)\in \hat{R}$. Furthermore, a subset of $A$ is an
$R$-submodule of $A$ if and only if it is a $\hat{R}$-submodule of $A$.

Next, we prove our main result in the special case that $R$ is local.

\begin{lemma}\label{26} Let $(R,\fm)$ be a local ring and $A$ an Artinian $R$-module of finite Gorenstein
injective dimension. Then $A$ admits a Gorenstein injective envelope $\G(A)$ which is Artinian and $\id_R(\G(A)/A)
<\infty$.
\end{lemma}

\begin{prf} First of all note that \cite[Lemma 3.6]{Sa} implies that $\Gid_{\hat{R}}A<\infty$. On the other
hand, for any Artinian ${\hat{R}}$-module $B$, it follows that $B$ is also Artinian as an $R$-module and $\Hom_R(B,B)=\Hom_{\hat{R}}(B,B)$. Hence, in view of Lemma \ref{24} and Lemma \ref{23} ii), we may and do
assume that $R$ is complete.

Denote the exact functor $\Hom_R(-,\E_R(R/\fm))$ by $(-)^{\vee}$. Then $N:=A^{\vee}$ is a finitely generated
$R$-module and $N^{\vee}\cong A$.  By \cite[Theorem 4.16]{CFoH} and \cite[Proposition 3.8 (a)]{CFrH}, one sees
that $\Gpd_RN=\Gid_RA<\infty.$ Hence, by Lemma \ref{21}, $N$ possesses a Gorenstein projective cover
$\phi:V\lo N$ such that $V$ is finitely generated and $K:=\Ker \phi$ has finite projective dimension. From
the exact sequence $$0\lo K\hookrightarrow V\overset{\phi}\lo N\lo 0,$$ we get the exact sequence $$0\lo A\overset{\lambda}\lo V^{\vee}\lo K^{\vee}\lo 0. \     \   \    \  (*)$$ Applying \cite[Theorem 4.16]{CFoH}
and \cite[Proposition 3.8 (a)]{CFrH} again implies that $V^{\vee}$ is an Artinian Gorenstein injective
$R$-module. We show that $\lambda$ is a Gorenstein injective envelope of $A$. As $\pd_RK<\infty$,  one
obtains $$\id_R(\Coker \lambda)=\id_RK^{\vee}<\infty.$$ For any Gorenstein injective $R$-module $G$,
\cite[Lemma 2.2]{CFrH} yields that $\Ext_R^i(K^{\vee},G)=0$ for all $i>0$, and hence from $(*)$, we
can deduce the exact sequence $$\Hom_R(V^{\vee},G)\lo \Hom_R(A,G)\lo \Ext_R^1(K^{\vee},G)=0.$$
Thus $\lambda:A\lo V^{\vee}$ is a Gorenstein injective preenvelope of $A$ and $\id_R(\Coker \lambda)
<\infty$.

Let $\tau:V\lo (V^{\vee})^{\vee}$ denote the natural $R$-isomorphism. Let $f:V^{\vee}\lo V^{\vee}$ be
an $R$-homomorphism such that $f\lambda=\lambda$. Then, we have $\phi(\tau^{-1}f^{\vee}\tau)=\phi$,
which yields that $\tau^{-1}f^{\vee}\tau$ is an automorphism. Therefore, $f$ is an automorphism, and so
$V^{\vee}$ is a Gorenstein injective envelope of $A$.
\end{prf}

The next is our last preliminary lemma for proving our main result.

\begin{lemma}\label{28} Let $S$ be a multiplicatively closed subset of $R$ and $L$ a $S^{-1}R$-module. Then $L$
is Gorenstein injective as an $R$-module if and only if it is Gorenstein injective as a $S^{-1}R$-module.
More generally, for any $S^{-1}R$-module $M$, we have $\Gid_RM=\Gid_{S^{-1}R}M$.
\end{lemma}

\begin{prf} First, assume that $L$ is Gorenstein injective as an $R$-module. Let $N$ be a $S^{-1}R$-
modules and $E^{\bullet}$ an injective resolution of $N$. For any $R$-module $M$, one has the following
natural isomorphisms of complexes:
$$\begin{array}{ll}\Hom_R(M,E^{\bullet})&\cong \Hom_R(M,\Hom_{S^{-1}R}(S^{-1}R,E^{\bullet}))\\
&\cong \Hom_{S^{-1}R}(M\otimes_RS^{-1}R,E^{\bullet})\\
&\cong \Hom_{S^{-1}R}(S^{-1}M,E^{\bullet}).
\end{array}$$
By bringing to mind that any injective $S^{-1}R$-module is also injective as an $R$-module, it follows
that $E^{\bullet}$ is also an injective resolution of $N$ as an $R$-module, and so $\Ext_R^i(M,N)\cong \Ext_{S^{-1}R}^i(S^{-1}M,N)$ for all $i\geq 0$. In particular, if $M$ is a $S^{-1}R$-module, then we have  $\Ext_R^i(M,N)\cong \Ext_{S^{-1}R}^i(M,N)$ for all $i\geq 0$. In what follows, we will use these
isomorphisms without any further comment. Let $$0\lo L\overset{d^{-1}}\lo I^0\overset {d^{0}}\lo I^1\lo
\cdots$$ be an augmented injective resolution of $L$ as a $S^{-1}R$-module and $I$ be an injective
$S^{-1}R$-module. Let $k\geq 0$ and $j>0$ be two integers. From the exact sequence $0\lo \im d^{k-1}\lo
I^k\lo \im d^k\lo 0,$ we deduce the $S^{-1}R$-isomorphisms $\Ext_{S^{-1}R}^{j}(I,\im d^k)\cong \Ext_{S^{-1}R}^{j+1}(I,\im d^{k-1}).$ Applying these isomorphisms successively for $k=i, i-1,\ldots, 0$
imply that
$$\begin{array}{ll}\Ext_{S^{-1}R}^{j}(I,\im d^i)&\cong \Ext_{S^{-1}R}^{j+i+1}(I,L)\\
&\cong \Ext_R^{j+i+1}(I,L)\\
&=0.$$
\end{array}$$
Recall by \cite[Theorem 3.3]{CFrH}, a module $E$ over a commutative Noetherian ring $T$ is Gorenstein
injective if and only if $\Gid_TE<\infty$ and $\Ext^i_T(\widehat{I},E)=0$ for all injective $T$-modules
$\widehat{I}$ and all $i>0$. By \cite[Theorem 5.4.1 and Corollary 7.2.3]{EJ}, the $S^{-1}R$-module $L$
admits a special injective cover $\phi_{0}:I_0\lo L$. Since $L$ is a Gorenstein injective $R$-module,
there are an injective $R$-module $\tilde{I}$ and an $R$-epimorphism $\tilde{I}\lo L$. Localizing at
$S$ yields a $S^{-1}R$-epimorphism $S^{-1}\tilde{I}\lo L.$ Hence, $\phi_{0}$ is an epimorphism. From
the exact sequence $0\lo \Ker \phi_{0}\hookrightarrow I_0\overset{\phi_{0}}\lo L\lo 0$, we deduce that
$\Gid_R(\Ker \phi_{0})<\infty$. On the other hand, as $\phi_{0}$ is special, one has $$\Ext_R^i(I',\Ker
\phi_{0})\cong \Ext_{S^{-1}R}^i(S^{-1}I',\Ker \phi_{0})=0$$ for all injective $R$-modules $I'$ and all
$i>0$. Thus, $\Ker \phi_{0}$ is a Gorenstein injective $R$-module. By repeating the above argument for
$\Ker \phi_{0}$ instead of $L$, we get a $S^{-1}R$-epimorphism $\phi_{1}:I_1\lo \Ker \phi_{0}$ where
$I_1$ is an injective $S^{-1}R$-module and and $\Ker \phi_{1}$ is a Gorenstein injective $R$-module.
Denote the composition $I_1\overset{\phi_{1}}\lo \Ker \phi_{0}\hookrightarrow I_0$ by $d_1$. Continuing
in this way, we find injective $S^{-1}R$-modules $I_0, I_1, I_2, \ldots $ and $S^{-1}R$-homomorphisms $I_n\overset{d_n}\lo I_{n-1}; \  \ n\geq 1$. Now, denote the sequence $$\cdots \lo I_n\overset{d_n}\lo
I_{n-1}\lo \cdots \lo I_0\overset{d^{-1}\phi_{0}}\lo I^0\overset{d^0}\lo I^1\lo \cdots$$ by $I^{\bullet}$.
Then $I^{\bullet}$ is an exact sequence of injective $S^{-1}R$-modules, $L\cong \im (I_0\lo I^0)$ and  $\Hom_{S^{-1}R}(I,I^{\bullet})$ is exact for all injective $S^{-1}R$-modules $I$. Therefore, $L$ is a
Gorenstein injective $S^{-1}R$-module.

Next, assume that $L$ is Gorenstein injective as an $S^{-1}R$-module. Then, there is an exact sequence
$$I^{\bullet}: \cdots \lo I_1\lo I_0\lo I^0\lo I^1\lo \cdots$$ of injective $S^{-1}R$-modules such
that $L\cong \im(I_0\lo I^0)$ and that $\Hom_{S^{-1}R}(I,I^{\bullet})$ is exact for all injective
$S^{-1}R$-modules $I$. Then $I^{\bullet}$ is also an exact sequence of injective $R$-modules. Let
$J$ be an injective $R$-module. We have the following natural isomorphisms of complexes: $$\Hom_R
(J,I^{\bullet})\cong \Hom_R(J,\Hom_{S^{-1}R}(S^{-1}R,I^{\bullet}))\cong \Hom_{S^{-1}R}(J\otimes_RS^{-1}R,I^{\bullet}).$$ But, $J\otimes_RS^{-1}R$ is an injective
$S^{-1}R$-module, and so the complex $\Hom_R(J,I^{\bullet})$ is exact. This shows that $L$ is a
Gorenstein injective $R$-module.

The last assertion follows immediately by Lemma \ref{23} i).
\end{prf}

We can record the following immediate corollary.

\begin{corollary}\label{29} Assume that every direct limit of Gorenstein injective $R$-modules is Gorenstein
injective. Let $S$ be a multiplicatively closed subset of $R$ and $M$ a Gorenstein injective
$R$-module. Then $S^{-1}M$ is a Gorenstein injective $S^{-1}R$-module.
\end{corollary}

\begin{prf} As $S^{-1}R$ is a flat $R$-module, the Govorov-Lazard Theorem implies that $S^{-1}R$ is
the direct limit of a direct system $\{F_\theta\}_{\theta\in \Theta}$ of finitely generated free
$R$-modules. One can easily check that $F_\theta\otimes_RM$ is a Gorenstein injective $R$-module for
all $\theta\in \Theta$. Thus, by our assumption, $S^{-1}M\cong \underset{\theta}{\varinjlim}
(F_\theta\otimes_RM)$ is a Gorenstein injective $R$-module. Next, by Lemma \ref{28}, we have $$\Gid_{S^{-1}R}S^{-1}M=\Gid_RS^{-1}M=0,$$ and so $S^{-1}M$ is is a Gorenstein injective
$S^{-1}R$-module, as required.
\end{prf}

Now, finally we are ready to prove our theorem.

\begin{theorem}\label{210} Let $A$ be an Artinian $R$-module of finite Gorenstein injective dimension.
Then $A$ admits a Gorenstein injective envelope $\G_R(A)$ which is Artinian and $\id_R(\G_R(A)/A)<\infty$.
\end{theorem}

\begin{prf} Let $\fm\in \Supp_RA$. Since $A$ is an Artinian $R$-module, there exists an exact sequence
$$0\lo A\lo E^0\lo E^1\lo \cdots \lo E^i\lo \cdots $$ in which each $E^i$ is an Artinian injective
$R$-module. Then, by \cite[Theorem 3.3]{CFrH}, it turns out that $\Ker(E^i\lo E^{i+1})$ is a Gorenstein
injective $R$-module for all $i\geq \Gid_RA$. Now, $$0\lo A_{\fm}\lo E^0_{\fm}\lo E^1_{\fm}\lo \cdots
\lo E^i_{\fm}\lo \cdots$$ is an injective resolution of the $R$-module $A_{\fm}$. By \cite[Lemma 2.1 (ii)]{D},
for each $i\geq 0$, $\Ker(E^i_{\fm}\lo E^{i+1}_{\fm})$ is a direct summand of $\Ker(E^i\lo E^{i+1})$. Hence, $\Ker(E^i_{\fm}\lo E^{i+1}_{\fm})$ is a Gorenstein injective $R$-module for
all $i\geq \Gid_RA$. Thus, $\Gid_RA_{\fm}\leq \Gid_RA$, and so $\Gid_{R_{\fm}}A_{\fm}$ is finite by
Lemma \ref{28}. Now, by Lemma \ref{26}, $\G_{R_{\fm}}(A_{\fm})$ exists and it is an Artinian
$R_{\fm}$-module with $\id_{R_{\fm}}(\G_{R_{\fm}}(A_{\fm})/A_{\fm})<\infty$. It is easy
to check that any Artinian $R_{\fm}$-module $B$ is also Artinian as an $R$-module and $\Hom_R(B,B)=\Hom_{R_{\fm}}(B,B)$.  Thus, $\G_{R_{\fm}}(A_{\fm})$ is an Artinian $R$-module and Lemma
\ref{28} and Lemma \ref{23} ii) yield that $\G_R(A_{\fm})=\G_{R_{\fm}}(A_{\fm})$ and  $\id_R(\G_R(A_{\fm})/A_{\fm})
<\infty$.

Suppose that $\Supp_RA=\{\fm_1,\ldots, \fm_n\}$. By \cite[Lemma 2.1 (ii)]{D}, one has a natural $R$-isomorphism
$$A\cong A_{\fm_1}\oplus\cdots\oplus A_{\fm_n}.$$ By the preceding paragraph for each $i$, the
$R$-module $A_{{\fm}_i}$ possesses an Artinian Gorenstein injective envelope $\phi_i:A_{{\fm}_i}
\lo \G_R(A_{{\fm}_i})$ such that $\id_R(\G_R(A_{{\fm}_i})/A_{{\fm}_i})<\infty$. Now, $\underset{i=1}
{\overset{n}\oplus} \G_R(A_{{\fm}_i})$ is an Artinian Gorenstein injective $R$-module and $$\underset{i}
\oplus \phi_i:\underset{i=1}{\overset{n}\oplus} A_{{\fm}_i}\lo \underset{i=1}{\overset{n}\oplus}
\G_R(A_{{\fm}_i})$$ is an $R$-monomorphism with $\Coker (\underset{i}\oplus \phi_i)\cong  \underset{i=1}
{\overset{n}\oplus} \G_R(A_{{\fm}_i})/A_{{\fm}_i}$. Hence, $\underset{i}\oplus \phi_i$ is a Gorenstein
injective preenvelope and $\id_R(\Coker (\underset{i}\oplus \phi_i))<\infty$. Let $\fm$ and $\fn$ be
two distinct maximal ideals of $R$ and $B$ and $C$ two $R$-modules such that $\Supp_RB\subseteq \{\fm\}$
and $\Supp_RC\subseteq \{\fn\}$. Then, one can easily check that $\Hom_R(B,C)=0$. Thus, we have $$\Hom_R(\underset{i=1}{\overset{n}\oplus}\G_R(A_{{\fm}_i}),\underset{i=1}{\overset{n}\oplus}
\G_R(A_{{\fm}_i}))\cong \underset{i=1}{\overset{n}\oplus} \Hom_R(\G_R(A_{{\fm}_i}),\G_R(A_{{\fm}_i})),$$
which yields that $\underset{i}\oplus \phi_i$ is our desired Gorenstein injective envelope.
\end{prf}

\begin{corollary}\label{211} Let $R$ be a Gorenstein ring and $A$ an Artinian $R$-module. Then $A$ admits a
Gorenstein injective envelope $\G_R(A)$ which is Artinian and $\id_R(\G_R(A)/A)<\infty$.
\end{corollary}

\begin{prf} Suppose that $\Supp_RA=\{\fm_1,\ldots, \fm_n\}$. By \cite[Lemma 2.1 (ii)]{D}, one has a natural
$R$-isomorphism $$A\cong A_{\fm_1}\oplus\cdots\oplus A_{\fm_n}.$$ Since over a Gorenstein local
ring, any module has finite Gorenstein injective dimension, the proof of Theorem \ref{210} yields that $$\Gid_RA=\max\{\Gid_RA_{\fm_i}|i=1,\ldots,n\}<\infty.$$ Hence, the claim follows by Theorem \ref{210}
\end{prf}

\begin{acknowledgement} The authors wish to thank Professor Edgar E. Enochs for his interest and
comments on this work.
\end{acknowledgement}

%%%%%%%%%%%%%%%%%%%%%%%%%%%%%%%%%%%%%%%%%%%%%%%%%%%%%

\end{document}